\documentclass[11pt]{amsart}
\usepackage{amssymb}
\usepackage{amsmath}
\usepackage[active]{srcltx}
\usepackage{t1enc}
\usepackage[latin2]{inputenc}
\usepackage{verbatim}
\usepackage{amsmath,amsfonts,amssymb,amsthm}
\usepackage[mathcal]{eucal}
\usepackage{enumerate}
\usepackage[centertags]{amsmath}
\usepackage{graphics}

\newtheorem{theorem}{Theorem}

\newtheorem{lemma}{Lemma}
\newtheorem{remark}{Remark}

\newtheorem{corollary}{Corollary}

\begin{document}
\author{L-E. Persson, G. Tephnadze and G. Tutberidze}
\title[Fejér means]{On the boundedness of subsequences of Vilenkin-Fejér
means on the martingale Hardy spaces }
\address {L-E. Persson, Department of Computer Science and Computational Engineering, UiT -The Arctic University of Norway, P.O. Box 385, N-8505, Narvik, Norway and Department of Mathematics and Computer Science, Karlstad University, Sweden}
\email{lars.e.persson@uit.no \ \ \ larserik.persson@kau.se}
\address {G. Tephnadze, The University of Georgia, School of Science and Technology, 77a Merab Kostava St, Tbilisi, 0128, Georgia.}
\email{g.tephnadze@ug.edu.ge}
\address{G.Tutberidze, The University of Georgia, School of science and technology, 77a Merab Kostava St, Tbilisi 0128, Georgia and Department of Computer Science and Computational Engineering, UiT -The Arctic University of Norway, P.O. Box 385, N-8505, Narvik, Norway.}
\email{giorgi.tutberidze1991@gmail.com}
\thanks{Second author was supported by grant of Shota Rustaveli National Science Foundation of Georgia, no. YS-18-043 and third author was supported by grant of Shota Rustaveli National Science Foundation of Georgia, no. PHDF-18-476.}
\date{}
\maketitle

\begin{abstract}
In this paper we characterize subsequences of Fejér means with respect to Vilenkin systems, which are bounded from the Hardy space $H_{p}$ to the Lebesgue space $L_{p},$ for all $0<p<1/2.$ The result is in a sense sharp.
\end{abstract}

\date{}



\section{Introduction}

\bigskip In the one-dimensional case the weak (1,1)-type inequality for the
maximal operator of Fejér means 
\begin{equation*}
\sigma ^{\ast }f:=\sup_{n\in \mathbb{N}}\left\vert \sigma _{n}f\right\vert
\end{equation*}%
can be found in Schipp \cite{Sc} for Walsh series and in Pál, Simon \cite{PS}
for bounded Vilenkin series. Here, as usual, the symbol $\sigma_n$ denotes the Fejér mean with respect to the Vilenkin system (and thus also called the Vilenkin-Fejér means, see Section 2). 

Fujji \cite{Fu} and Simon \cite{Si2} verified
that $\sigma ^{\ast }$ is bounded from $H_{1}$ to $L_{1}$. Weisz \cite{We2}
generalized this result and proved boundedness of $\sigma ^{\ast }$ from the
martingale space $H_{p}$ to the Lebesgue space $L_{p}$ for $p>1/2$. Simon 
\cite{Si1} gave a counterexample, which shows that boundedness does not hold
for $0<p<1/2.$ A counterexample for $p=1/2$ was given by Goginava \cite%
{GoAMH} (see also \cite{BGG} and \cite{BGG2}). Weisz \cite{We4} proved that
the maximal operator of the Fejér means $\sigma ^{\ast }$ is bounded from
the Hardy space $H_{1/2}$ to the space $weak-L_{1/2}$. The boundedness of
weighted maximal operators are considered in \cite{GNCz}, 
\cite{tep2} and \cite{tep3}. 

Weisz \cite{We3} (see also \cite{We1}) also proved
that the following theorem is true:

\textbf{Theorem W:} \textbf{(Weisz)} Let $p>0.$ Then the maximal operator 
\begin{equation} \label{we}
\sigma ^{\nabla ,\ast }f=\underset{n\in \mathbb{N}}{\sup }\left\vert \sigma
_{M_{n}}f\right\vert
\end{equation}
where 
$
M_{0}:=1,\ M_{n+1}:=m_{n}M_{n}\  (n\in \mathbb{N})
$
and $m:=(m_{0},m_{1},\dots )$ be a sequences of
the positive integers not less than 2, which generate Vilenkin systems,
is bounded from the Hardy space $H_{p}$ to the space $L_{p}.$

In \cite{PT1} the result of Weisz was generalized and it was found the maximal subspace $S\subset \mathbb{N}$ of positive numbers, for which the restricted maximal operator on this subspace
$
\underset{n\in S\subset\mathbb{N}}{\sup }\left\vert \sigma
_{n}f\right\vert
$
of Fejér means  is bounded from the Hardy space $H_{p}$ to the space $L_{p}$ for all $0<p\leq 1/2.$ The new theorem (Theorem 1) in this paper show in particular that this result is in a sense sharp. In particular, for every natural number $n=\sum_{k=0}^{\infty }n_{k}M_{k},$ where $n_{k}\in Z_{m_{k}}\ (k\in \mathbb{N%
}_{+})$ we define numbers 
\begin{equation*}
\left\langle n\right\rangle :=\min \{j\in \mathbb{N}:n_{j}\neq 0\},
 \ \ \ \ \ \left\vert n\right\vert :=\max \{j\in \mathbb{N}:n_{j}\neq 0\},  \ \ \ \ \rho \left( n\right) =\left\vert n\right\vert -\left\langle
n\right\rangle
\end{equation*}
and prove that 
$$ S=\{n\in \mathbb{N}:\rho \left( n\right)\leq c<\infty. \}$$

Since $\rho (M_n)=0$ for all $n\in \mathbb{N}$ we obtain that 
$\{M_n: n\in \mathbb{N}\}\subset S$ and that follows i.e. that result of Weisz \cite{We3} (see also \cite{We1}) that restricted maximal operator \eqref{we} is bounded from the Hardy space $H_{p}$ to the space $L_{p}.$

The main aim of this paper is to generalize Theorem W and find the maximal
subspace of positive numbers, for which the restricted maximal operator of
Fejér means in this subspace is bounded from the Hardy space $H_{p}$ to the
space $L_{p}$ for all $0<p\leq 1/2.$ As applications, both some well-known
and new results are pointed out.

This paper is organized as follows: In order not to disturb our discussions
later on some preliminaries (definitions, notations and lemmas) are presented in Section 2. The main result (Theorem 1) and some of its consequences can be found in Section 3. The detailed proof of Theorem 1 is given in Section 4.

\section{Preliminaries}

Denote by $\mathbb{N}_{+}$ the set of the positive integers, $\mathbb{N}:=%
\mathbb{N}_{+}\cup \{0\}.$ Let $m:=(m_{0},m_{1},\dots )$ be a sequence of
the positive integers not less than 2. Denote by $Z_{m_{n}}:=\{0,1,\ldots
,m_{n}-1\}$ the additive group of integers modulo $m_{n}$. Define the group $%
G_{m}$ as the complete direct product of the groups $Z_{m_{n}}$ with the
product of the discrete topologies of $Z_{m_{n}}`$s. In this paper we
discuss bounded Vilenkin groups, i.e. the case when $\sup_{n\in \mathbb{N}%
}m_{n}<\infty .$

The direct product $\mu $ of the measures $\mu _{n}\left( \{j\}\right)
:=1/m_{n},\ (j\in Z_{m_{n}})$ is the Haar measure on $G_{m}$ with $\mu
\left( G_{m}\right) =1.$

The elements of $G_{m}$ are represented by sequences 
\begin{equation*}
x:=\left( x_{0},x_{1},\ldots,x_{n},\ldots \right),\ \left( x_{n}\in
Z_{m_{n}}\right).
\end{equation*}

It is easy to give a base for the neighbourhood of $G_{m}:$

\begin{equation*}
I_{0}\left( x\right) :=G_{m},\ I_{n}(x):=\{y\in G_{m}\mid
y_{0}=x_{0},\ldots,y_{n-1}=x_{n-1}\}\,\,\left( x\in G_{m},\text{ }n\in 
\mathbb{N}\right).
\end{equation*}

Set $I_{n}:=I_{n}\left( 0\right) ,$ for $n\in \mathbb{N}_{+}$ \ and

\begin{equation*}
e_{n}:=\left( 0,\dots ,0,x_{n}=1,0,\dots \right) \in G_{m}\qquad \left( n\in 
\mathbb{N}\right) .
\end{equation*}

Denote 
\begin{equation*}
I_{N}^{k,l}:=\left\{ 
\begin{array}{l}
\text{ }I_{N}(0,\dots ,0,x_{k}\neq 0,0,\dots ,0,x_{l}\neq 0,x_{l+1,\dots ,%
\text{ }}x_{N-1\text{ }}),\text{ \ \ }k<l<N, \\ 
\text{ }I_{N}(0,\dots ,0,x_{k}\neq 0,0,\dots ,0),\text{ \qquad }l=N.%
\end{array}%
\text{ }\right.
\end{equation*}

It is easy to show that 
\begin{equation}
\overline{I_{N}}=\left( \overset{N-2}{\underset{i=0}{\bigcup }}\overset{N-1}{%
\underset{j=i+1}{\bigcup }}I_{N}^{i,j}\right) \bigcup \left( \underset{i=0}{%
\bigcup\limits^{N-1}}I_{N}^{i,N}\right), \ \ \ n=2,3,...  \label{2}
\end{equation}

\bigskip If we define the so-called generalized number system based on $m$
in the following way : 
\begin{equation*}
M_{0}:=1,\ M_{n+1}:=m_{n}M_{n}\ \ \ (n\in \mathbb{N}),
\end{equation*}%
then every $n\in \mathbb{N}$ can be uniquely expressed as $%
n=\sum_{k=0}^{\infty }n_{k}M_{k},$ where $n_{k}\in Z_{m_{k}}\ (k\in \mathbb{N%
}_{+})$ and only a finite number of $n_{k}`$s differ from zero. Let 
\begin{equation*}
\left\langle n\right\rangle :=\min \{j\in \mathbb{N}:n_{j}\neq 0\}\text{ \ \
and \ \ \ }\left\vert n\right\vert :=\max \{j\in \mathbb{N}:n_{j}\neq 0\},
\end{equation*}%
that is $M_{\left\vert n\right\vert }\leq n\leq M_{\left\vert n\right\vert
+1}.$ Set $\rho \left( n\right) =\left\vert n\right\vert -\left\langle
n\right\rangle ,$ \ for \ all \ \ $n\in \mathbb{N}.$

Next, we introduce on $G_{m}$ an orthonormal system, which is called the
Vilenkin system. At first, we define the complex-valued function $%
r_{k}\left( x\right) :G_{m}\rightarrow \mathbb{C},$ the generalized
Rademacher functions, by 
\begin{equation*}
r_{k}\left( x\right) :=\exp \left( 2\pi ix_{k}/m_{k}\right) ,\text{ }\left(
i^{2}=-1,x\in G_{m},\text{ }k\in \mathbb{N}\right) .
\end{equation*}

Now, define the Vilenkin system $\,\,\,\psi :=(\psi _{n}:n\in\mathbb{N})$ on 
$G_{m}$ as: 
\begin{equation*}
\psi _{n}(x):=\prod\limits_{k=0}^{\infty }r_{k}^{n_{k}}\left( x\right) \ \ \
\left( n\in\mathbb{N}\right).
\end{equation*}

Specifically, we call this system the Walsh-Paley system, when $m\equiv 2.$

The norms (or quasi-norms) of the spaces $L_{p}(G_{m})$ and $%
weak-L_{p}\left( G_{m}\right) $ $\ \left( 0<p<\infty \right) $ are
respectively defined by 
\begin{equation*}
\left\Vert f\right\Vert _{p}^{p}:=\int_{G_{m}}\left\vert f\right\vert
^{p}d\mu ,\text{ \ \ }\left\Vert f\right\Vert _{L_{p,\infty}}^{p}:=\underset{%
\lambda >0}{\sup }\,\lambda ^{p}\mu \left( f>\lambda \right) <\infty .
\end{equation*}%
\qquad

The Vilenkin system is orthonormal and complete in $L_{2}\left( G_{m}\right) 
$ (see \cite{Vi}).

If $\ f\in L_{1}\left( G_{m}\right) $ we can define Fourier coefficients,
partial sums, Dirichlet kernels, Fejér means, Fejér kernels with respect to
the Vilenkin system in the usual manner: 
\begin{equation*}
\widehat{f}\left( k\right) :=\int_{G_{m}}f\overline{\psi }_{k}d\mu \text{%
\thinspace\ \ \ }\left( \text{ }k\in \mathbb{N}\right) ,
\end{equation*}%
\begin{eqnarray*}
S_{n}f &:&=\sum_{k=0}^{n-1}\widehat{f}\left( k\right) \psi _{k},\ \text{%
\qquad }D_{n}:=\sum_{k=0}^{n-1}\psi _{k\text{ }}\text{ \ \ \ \ \ \ }\left( 
\text{ }n\in \mathbb{N}_{+}\text{ }\right) , \\
\sigma _{n}f &:&=\frac{1}{n}\sum_{k=0}^{n-1}S_{k}f,\text{\qquad\ \ \ }K_{n}:=%
\frac{1}{n}\overset{n-1}{\underset{k=0}{\sum }}D_{k}\text{ \ \ }\left( \text{
}n\in \mathbb{N}_{+}\text{ }\right) .
\end{eqnarray*}

Recall that (see e.g. \cite{AVD}) 
\begin{equation}
\quad \hspace*{0in}D_{M_{n}}\left( x\right) =\left\{ 
\begin{array}{l}
\text{ }M_{n},\text{ \thinspace \thinspace if \thinspace \thinspace }x\in
I_{n}, \\ 
\text{ }0,\text{ \ \ \thinspace \thinspace \thinspace if \ \thinspace
\thinspace }x\notin I_{n},%
\end{array}%
\right.  \label{3}
\end{equation}%
and \ 
\begin{equation}
D_{s_{n}M_{n}}=D_{s_{n}M_{n}}\sum_{k=0}^{s_{n}-1}\psi
_{kM_{n}}=D_{M_{n}}\sum_{k=0}^{s_{n}-1}r_{n}^{k},  \label{9dn}
\end{equation}%
where $n\in \mathbb{N}$ and $1\leq s_{n}\leq m_{n}-1.$

The $\sigma $-algebra generated by the intervals $\left\{ I_{n}\left(
x\right) :x\in G_{m}\right\} $ will be denoted by $\digamma _{n}$ $\left(
n\in \mathbb{N}\right) .$ Denote by $f=\left( f^{\left( n\right) },n\in 
\mathbb{N}\right) $ a martingale with respect to $\digamma _{n}$ $\left(
n\in \mathbb{N}\right) $ (for details see e.g. \cite{We1}). The maximal
function of a martingale $f$ is defined by \qquad 
\begin{equation*}
f^{\ast }=\sup_{n\in \mathbb{N}}\left\vert f^{\left( n\right) }\right\vert .
\end{equation*}

In the case $f\in L_{1}(G_{m}),$ the maximal functions are just also given by 
\begin{equation*}
f^{\ast }\left( x\right) =\sup_{n\in \mathbb{N}}\frac{1}{\left\vert
I_{n}\left( x\right) \right\vert }\left\vert \int_{I_{n}\left( x\right)
}f\left( u\right) \mu \left( u\right) \right\vert .
\end{equation*}

For $0<p<\infty $ the Hardy martingale spaces $H_{p}\left( G_{m}\right) $
consist of all martingales $f$, for which 
\begin{equation*}
\left\Vert f\right\Vert _{H_{p}}:=\left\Vert f^{\ast }\right\Vert
_{p}<\infty .
\end{equation*}

If $f\in L_{1}(G_{m}),$ then it is easy to show that the sequence $\left(
S_{M_{n}}\left( f\right) :n\in \mathbb{N}\right) $ is a martingale. If $%
f=\left( f^{\left( n\right) },n\in \mathbb{N}\right) $ is a martingale, then
the Vilenkin-Fourier coefficients must be defined in a slightly different
manner: $\qquad \qquad $ 
\begin{equation*}
\widehat{f}\left( i\right) :=\lim_{k\rightarrow \infty
}\int_{G_{m}}f^{\left( k\right) }\left( x\right) \overline{\psi }_{i}\left(
x\right) d\mu \left( x\right) .
\end{equation*}

The Vilenkin-Fourier coefficients of $f\in L_{1}\left( G_{m}\right) $ are
the same as those of the martingale $\left( S_{M_{n}} f :n\in 
\mathbb{N}\right) $ obtained from $f$.

A bounded measurable function $a$ is said to be a p-atom if there exists an
interval $I$, such that%
\begin{equation*}
\int_{I}ad\mu =0,\text{ \ \ }\left\Vert a\right\Vert _{\infty }\leq \mu
\left( I\right) ^{-1/p},\text{ \ \ supp}\left( a\right) \subset I.\qquad
\end{equation*}%
\qquad

For the proof of the main result (Theorem 1) we need the following Lemmas:

\begin{lemma}[see e.g. \protect\cite{We3}]
	\label{lemma1} A martingale $f=\left( f^{\left( n\right) },n\in \mathbb{N}%
	\right) $ is in $H_{p}\left( 0<p\leq 1\right) $ if and only if there exist a
	sequence $\left( a_{k},k\in \mathbb{N}\right) $ of p-atoms and a sequence $%
	\left( \mu _{k},k\in \mathbb{N}\right) $ of real numbers such that for every 
	$n\in \mathbb{N}:$%
	\begin{equation}
	\qquad \sum_{k=0}^{\infty }\mu _{k}S_{M_{n}}a_{k}=f^{\left( n\right) }
	\label{1}
	\end{equation}%
	and%
	\begin{equation*}
	\qquad \sum_{k=0}^{\infty }\left\vert \mu _{k}\right\vert ^{p}<\infty .
	\end{equation*}%
	Moreover, $\left\Vert f\right\Vert _{H_{p}}\backsim \inf \left(
	\sum_{k=0}^{\infty }\left\vert \mu _{k}\right\vert ^{p}\right) ^{1/p},$
	where the infimum is taken over all decomposition of $f$ of the form (\ref{1}%
	).
\end{lemma}

\begin{lemma}[see e.g. \protect\cite{We3}]
	\label{lemma2} Suppose that an operator $T$ is $\sigma $-linear and for some 
	$0<p\leq 1$%
	\begin{equation*}
	\int\limits_{\overset{-}{I}}\left\vert Ta\right\vert ^{p}d\mu \leq
	c_{p}<\infty ,
	\end{equation*}%
	for every $p$-atom $a$, where $I$ denotes the support of the atom. If $T$ is
	bounded from $L_{\infty \text{ }}$ to $L_{\infty },$ then 
	\begin{equation*}
	\left\Vert Tf\right\Vert _{p}\leq c_{p}\left\Vert f\right\Vert _{H_{p}}.
	\end{equation*}
\end{lemma}

\begin{lemma}[see \protect\cite{gat}]
	\label{lemma3} Let $n>t,$ $t,n\in \mathbb{N},$ $x\in I_{t}\backslash $ $%
	I_{t+1}$. Then 
	\begin{equation*}
	K_{M_{n}}\left( x\right) =\left\{ 
	\begin{array}{ll}
	0, & \text{if }x-x_{t}e_{t}\notin I_{n}, \\ 
	\frac{M_{t}}{1-r_{t}\left( x\right) }, & \text{if }x-x_{t}e_{t}\in I_{n}.%
	\end{array}%
	\right.
	\end{equation*}
\end{lemma}

\begin{lemma}[see \protect\cite{tep3}]
	\label{lemma5b}\ Let $x\in I_{N}^{i,j},$ $i=0,\dots ,N-1,$ $j=i+1,\dots ,N$.
	Then%
	\begin{equation*}
	\int_{I_{N}}\left\vert K_{n}\left( x-t\right) \right\vert d\mu \left(
	t\right) \leq \frac{cM_{i}M_{j}}{M_{N}^{2}},\text{ \ \ \ for \ }n\geq M_{N}.
	\end{equation*}
\end{lemma}

\begin{lemma} [see \cite{PT1}]
	\label{lemma8}Let $n\in \mathbb{N}.$ Then%
	\begin{equation}
	\left\vert K_{n}\left( x\right) \right\vert \leq \frac{c}{n}%
	\sum_{l=\left\langle n\right\rangle }^{\left\vert n\right\vert
	}M_{l}\left\vert K_{M_{l}}\right\vert \leq c\sum_{l=\left\langle
		n\right\rangle }^{\left\vert n\right\vert }\left\vert K_{M_{l}}\right\vert
	\label{10k}
	\end{equation}%
	and 
	\begin{equation}
	\left\vert nK_{n}\right\vert \geq \frac{M_{\left\langle n\right\rangle }^{2}%
	}{2\pi \lambda },\text{ \ \ \ }x\in I_{\left\langle n\right\rangle +1}\left(
	e_{\left\langle n\right\rangle -1}+e_{\left\langle n\right\rangle }\right) ,
	\label{9k}
	\end{equation}%
	where $\lambda :=\sup m_{n}.$
\end{lemma}

\section{The Main Result and applications}

Our main result reads:

\begin{theorem}
\label{theorem1}a) Let $0<p<1/2,$ $f\in H_{p}.$ Then there exists an
absolute constant $c_{p}$, depending only on $p$, such that 
\begin{equation*}
\text{ }\left\Vert \sigma _{n_{k}}f\right\Vert _{H_{p}}\leq \frac{%
c_{p}M_{\left\vert n_{k}\right\vert }^{1/p-2}}{M_{\left\langle
n_{k}\right\rangle }^{1/p-2}}\left\Vert f\right\Vert _{H_{p}}.
\end{equation*}

\textit{b) (sharpness) Let }$0<p<1/2$ \textit{and} $\Phi \left( n\right) $ \textit{be
any nondecreasing function,\ such that } 
\begin{equation} \label{31aaa}
\sup_{k\in \mathbb{N}}\rho\left( n_{k}\right) =\infty ,\text{ \ \ }\overline{\underset{%
k\rightarrow \infty }{\lim }}\frac{M_{\left\vert n_{k}\right\vert }^{1/p-2}}{%
M_{\left\langle n_{k}\right\rangle}^{1/p-2}\Phi\left( n_{k}\right) }
=\infty.  
\end{equation}%
\textit{Then there exists a martingale }$f\in H_{p},$\textit{\ such that} 
\begin{equation*}
\underset{k\in \mathbb{N}}{\sup }\left\Vert \frac{\sigma _{n_{k}}f}{\Phi \left(
n_{k}\right) }\right\Vert _{L_{p,\infty }}=\infty .
\end{equation*}
\end{theorem}

\begin{corollary}
	\label{corollary1} Let $0<p<1/2,$ and  $f\in H_{p}.$ Then there exists
	an absolute constant $c_{p}$, depending only on $p$, such that 
	\begin{equation*}
	\text{ }\left\Vert \sigma _{n_k}f\right\Vert _{H_{p}}\leq c_{p}\left\Vert
	f\right\Vert _{H_{p}},\text{ \ \ }k\in \mathbb{N}
	\end{equation*}
	if and only if 
	\begin{equation*}
	\sup_{k\in \mathbb{N}}\rho\left( n_{k}\right)<c<\infty.
	\end{equation*}
\end{corollary}

As an application we also obtain the previous mentioned result by Weisz \cite
{We1}, \cite{We3} (Theorem W).

\begin{corollary}\label{corollary2} 
	Let $0<p<1/2,$ \  $f\in H_{p}.$ Then there exists
an absolute constant $c_{p}$, depending only on $p$, such that 
\begin{equation*}
\text{ }\left\Vert \sigma _{M_{n}}f\right\Vert _{H_{p}}\leq c_{p}\left\Vert
f\right\Vert _{H_{p}},\text{ \ \ }n\in \mathbb{N}.
\end{equation*}
\end{corollary}
On the other hand, the following unexpected result is true:
\begin{corollary}
\label{corollary3} \textit{a)} Let $0<p<1/2,$ \  $f\in H_{p}.$ Then there exists an absolute constant $c_{p}$, depending only on $p$, such that 
\begin{equation*}
\text{ }\left\Vert \sigma _{M_{n}+1}f\right\Vert _{H_{p}}\leq c_{p}M^{1/p-2}_{n}\left\Vert
f\right\Vert _{H_{p}},\text{ \ \ }n\in \mathbb{N}.
\end{equation*}
\textit{b)} Let $0<p<1/2$ \textit{and} $\Phi \left( n\right) $ \textit{be any nondecreasing function, such that} 
\begin{equation*} 
\overline{\underset{	k\rightarrow \infty }{\lim }}\frac{M_{k }^{1/p-2}}{\Phi\left( k\right)}
=\infty.  
\end{equation*}%
\textit{Then there exists a martingale }$f\in H_{p},$\textit{\ such that} 
\begin{equation*}
\underset{k\in \mathbb{N}}{\sup }\left\Vert \frac{\sigma _{M_{k}+1}f}{\Phi \left(k\right) }\right\Vert _{L_{p,\infty }}=\infty .
\end{equation*}
\end{corollary}
\begin{remark}
From Corollary \ref{corollary2} we obtain that $ \sigma _{M_{n}} $ are bounded from $ H_p $ to $ H_p $, but from Corollary \ref{corollary3} we conclude that $ \sigma _{M_{n}+1} $ are not bounded from $ H_p $ to $ H_p $. The main reason is that Fourier coefficients of martingales $ f\in H_p $ are not uniformly bounded (for details see e.g. \cite{tep8}).
\end{remark}
In the next corollary we state some estimates for the Walsh system only to clearly see the difference of divergence rates for the various subsequences:
\begin{corollary}
\label{corollary4} \textit{a)}Let $0<p<1/2,$ \  $f\in H_{p}.$ Then there exists an absolute constant $c_{p}$, depending only on $p$, such that 
\begin{equation} \label{aaa}
\text{ }\left\Vert \sigma _{2^{n}+1}f\right\Vert _{H_{p}}\leq c_{p}2^{(1/p-2)n}\left\Vert
f\right\Vert _{H_{p}},\text{ \ \ }n\in \mathbb{N}
\end{equation}
and
\begin{equation} \label{aaa1}
\text{ }\left\Vert \sigma _{2^{n}+1}f\right\Vert _{H_{p}}\leq c_{p}2^{\frac{(1/p-2)n}{2}}\left\Vert
f\right\Vert _{H_{p}},\text{ \ \ }n\in \mathbb{N}.
\end{equation}
\textit{b)} The rates $2^{(1/p-2)n} $ and  $ 2^{\frac{(1/p-2)n}{2}}$ in inequalities (\ref{aaa}) and (\ref{aaa1}) are sharp in the same sense as in Theorem \ref{theorem1}.
\end{corollary}

\section{Proof of Theorem 1}

\begin{proof}
a) Since 
\begin{equation}
\sup_{n\in \mathbb{N}}\int_{G_{m}}\left\vert K_{n}\left( x\right)
\right\vert d\mu \left( x\right) \leq c<\infty,   \label{4}
\end{equation}
we obtain that 
$$\frac{M_{\left\langle n_{k}\right\rangle
}^{1/p-2}\left\vert \sigma _{n_{k}}a\left( x\right) \right\vert }{%
M_{\left\vert n_{k}\right\vert }^{1/p-2}}$$   
is bounded from $L_{\infty}$ to $L_{\infty}.$ According to Lemma \ref{lemma2} we find that the proof of
Theorem \ref{theorem1} will be complete, if we show that%
\begin{equation*}
\int_{\overline{I_{N}}}\left\vert \frac{M_{\left\langle n_{k}\right\rangle
}^{1/p-2}\sigma _{n_{k}}a\left( x\right) }{M_{\left\vert n_{k}\right\vert
}^{1/p-2}}\right\vert ^{p}<c<\infty ,
\end{equation*}%
for every $p$-atom $a,$ with support$\ I$ and $\mu \left( I\right)
=M_{N}^{-1}.$ We may assume that $I=I_{N}.$ It is easy to see that $\sigma
_{n_{k}}\left( a\right) =0$ when $n_{k}\leq M_{N}.$ Therefore, we can
suppose that $n_{k}>M_{N}$.

Since $\left\Vert a\right\Vert _{\infty }\leq M_{N}^{1/p}$ we find that%
\begin{eqnarray} \label{400}
&&\frac{M_{\left\langle n_{k}\right\rangle }^{1/p-2}\left\vert \sigma
_{n_{k}}a\left( x\right) \right\vert }{M_{\left\vert n_{k}\right\vert
}^{1/p-2}}\leq \frac{M_{\left\langle n_{k}\right\rangle }^{1/p-2}}{%
M_{\left\vert n_{k}\right\vert }^{1/p-2}}\int_{I_{N}}\left\vert a\left(
t\right) \right\vert \left\vert K_{n_{k}}\left( x-t\right) \right\vert d\mu
\left( t\right)   \label{400a} \\
&\leq &\frac{M_{\left\langle n_{k}\right\rangle }^{1/p-2}\left\Vert
a\right\Vert _{\infty }}{M_{\left\vert n_{k}\right\vert }^{1/p-2}}%
\int_{I_{N}}\left\vert K_{n_{k}}\left( x-t\right) \right\vert d\mu \left(
t\right)   \notag \\
&\leq &\frac{M_{\left\langle n_{k}\right\rangle }^{1/p-2}M_{N}^{1/p}}{%
M_{\left\vert n_{k}\right\vert }^{1/p-2}}\int_{I_{N}}\left\vert
K_{n_{k}}\left( x-t\right) \right\vert d\mu \left( t\right)   \notag \\
&\leq &M_{\left\langle n_{k}\right\rangle }^{1/p-2}M_{\left\vert
n_{k}\right\vert }^{2}\int_{I_{N}}\left\vert K_{n_{k}}\left( x-t\right)
\right\vert d\mu \left( t\right) .  \notag
\end{eqnarray}

Without loss the generality we may assume that $i<j$. Let $x\in I_{N}^{i,j}\ 
$and $j<\left\langle n_{k}\right\rangle .$ Then $x-t\in I_{N}^{i,j}$ for $%
t\in I_{N}$ and, according to Lemma \ref{lemma3}, we obtain that 
\begin{equation*}
\left\vert K_{M_{l}}\left( x-t\right) \right\vert =0,\text{ \ for all }
\left\langle n_{k}\right\rangle \leq \text{ }l\leq \left\vert
n_{k}\right\vert.
\end{equation*}

By applying (\ref{400}) and (\ref{10k}) in Lemma \ref{lemma8}, for $ x \in I_{N}^{i,j},\text{ \ }0\leq i<j<\left\langle
n_{k}\right\rangle $ we get that 
\begin{eqnarray} \label{401} 
&&\frac{M_{\left\langle n_{k}\right\rangle }^{1/p-2}\left\vert \sigma
_{n_{k}}a\left( x\right) \right\vert }{M_{\left\vert n_{k}\right\vert
}^{1/p-2}} 
\leq M_{\left\langle n_{k}\right\rangle }^{1/p-2}M_{\left\vert
n_{k}\right\vert }^{2}\overset{\left\vert n_{k}\right\vert }{\underset{%
l=\left\langle n_{k}\right\rangle }{\sum }}\int_{I_{N}}\left\vert
K_{M_{l}}\left( x-t\right) \right\vert d\mu \left( t\right) =0.\text{ }
\end{eqnarray}

Let $x\in I_{N}^{i,j},\,$where $\left\langle n_{k}\right\rangle \leq j\leq N.
$ Then, in the view of Lemma \ref{lemma5b}, we have that 
\begin{equation*}
\int_{I_{N}}\left\vert K_{n_{k}}\left( x-t\right) \right\vert d\mu \left(
t\right) \leq \frac{cM_{i}M_{j}}{M_{N}^{2}}.
\end{equation*}

By using again (\ref{400}) we find that 
\begin{eqnarray}  \label{402}
\frac{M_{\left\langle n_{k}\right\rangle }^{1/p-2}\left\vert \sigma
_{n_{k}}a\left( x\right) \right\vert }{M_{\left\vert n_{k}\right\vert
}^{1/p-2}} &\leq &\frac{M_{\left\langle n_{k}\right\rangle
}^{1/p-2}M_{N}^{1/p}}{M_{\left\vert n_{k}\right\vert }^{1/p-2}}%
\int_{I_{N}}\left\vert K_{n_{k}}\left( x-t\right) \right\vert d\mu \left(
t\right)   \\
&\leq &\frac{M_{\left\langle n_{k}\right\rangle }^{1/p-2}M_{N}^{1/p}}{%
M_{\left\vert n_{k}\right\vert }^{1/p-2}}\frac{M_{i}M_{j}}{M_{N}^{2}}\leq
M_{\left\langle n_{k}\right\rangle }^{1/p-2}M_{i}M_{j}.  \notag
\end{eqnarray}

By combining (\ref{2}) and (\ref{400})-(\ref{402}) we get that 
\begin{eqnarray*}
&&\int_{\overline{I_{N}}}\left\vert \frac{M_{\left\langle n_{k}\right\rangle
}^{1/p-2}\left\vert \sigma _{n_{k}}a\left( x\right) \right\vert }{%
M_{\left\vert n_{k}\right\vert }^{1/p-2}}\right\vert ^{p}d\mu  \\
&=&\overset{N-2}{\underset{i=0}{\sum }}\overset{N-1}{\underset{j=i+1}{\sum }}%
\int_{I_{N}^{i,j}}\left\vert \frac{M_{\left\langle n_{k}\right\rangle
}^{1/p-2}\left\vert \sigma _{n_{k}}a\left( x\right) \right\vert }{%
M_{\left\vert n_{k}\right\vert }^{1/p-2}}\right\vert ^{p}d\mu \\ &+&\overset{N-1}{%
\underset{i=0}{\sum }}\int_{I_{N}^{k,N}}\left\vert \frac{M_{\left\langle
n_{k}\right\rangle }^{1/p-2}\left\vert \sigma _{n_{k}}a\left( x\right)
\right\vert }{M_{\left\vert n_{k}\right\vert }^{1/p-2}}\right\vert ^{p}d\mu 
\\
&\leq &\overset{\left\langle n_{k}\right\rangle -1}{\underset{i=0}{\sum }}%
\overset{N-1}{\underset{j=\left\langle n_{k}\right\rangle }{\sum }}%
\int_{I_{N}^{i,j}}\left\vert \frac{M_{\left\langle n_{k}\right\rangle
}^{1/p-2}\left\vert \sigma _{n_{k}}a\left( x\right) \right\vert }{%
M_{\left\vert n_{k}\right\vert }^{1/p-2}}\right\vert ^{p}d\mu  \\
&&+\overset{N-2}{\underset{i=\left\langle n_{k}\right\rangle }{\sum }}%
\overset{N-1}{\underset{j=i+1}{\sum }}\int_{I_{N}^{i,j}}\left\vert \frac{%
M_{\left\langle n_{k}\right\rangle }^{1/p-2}\left\vert \sigma
_{n_{k}}a\left( x\right) \right\vert }{M_{\left\vert n_{k}\right\vert
}^{1/p-2}}\right\vert ^{p}d\mu  \\
&&+\overset{N-1}{\underset{i=0}{\sum }}\int_{I_{N}^{i,N}}\left\vert \frac{%
M_{\left\langle n_{k}\right\rangle }^{1/p-2}\left\vert \sigma
_{n_{k}}a\left( x\right) \right\vert }{M_{\left\vert n_{k}\right\vert
}^{1/p-2}}\right\vert ^{p}d\mu  
\end{eqnarray*}
\begin{eqnarray*}
&\leq &\overset{\left\langle n_{k}\right\rangle -1}{\underset{i=0}{\sum }}%
\overset{N-1}{\underset{j=\left\langle n_{k}\right\rangle }{\sum }}%
\int_{I_{N}^{i,j}}\left\vert M_{\left\langle n_{k}\right\rangle
}^{1/p-2}M_{i}M_{j}\right\vert ^{p}d\mu +\overset{N-2}{\underset{%
i=\left\langle n_{k}\right\rangle }{\sum }}\overset{N-1}{\underset{j=i+1}{%
\sum }}\int_{I_{N}^{i,j}}\left\vert M_{\left\langle n_{k}\right\rangle
}^{1/p-2}M_{i}M_{j}\right\vert ^{p}d\mu  \\
&&+\overset{N-1}{\underset{i=0}{\sum }}\int_{I_{N}^{i,N}}\left\vert
M_{\left\langle n_{k}\right\rangle }^{1/p-2}M_{i}M_{N}\right\vert ^{p}d\mu 
\\
&\leq &c_{p}M_{\left\langle n_{k}\right\rangle }^{1-2p}\overset{\left\langle
n_{k}\right\rangle -1}{\underset{i=0}{\sum }}\overset{N-1}{\underset{%
j=\left\langle n_{k}\right\rangle }{\sum }}\frac{\left( M_{i}M_{j}\right)
^{p}}{M_{j}}+c_{p}M_{\left\langle n_{k}\right\rangle }^{1-2p}\overset{N-2}{%
\underset{i=\left\langle n_{k}\right\rangle }{\sum }}\overset{N-1}{\underset{%
j=i+1}{\sum }}\frac{\left( M_{i}M_{j}\right) ^{p}}{M_{j}} \\
&&+c_{p}M_{\left\langle n_{k}\right\rangle }^{1-2p}\underset{i=0}{\sum }%
\frac{\left( M_{i}M_{N}\right) ^{p}}{M_{N}} \\
&\leq &c_{p}M_{\left\langle n_{k}\right\rangle }^{1-2p}\overset{\left\langle
n_{k}\right\rangle }{\underset{i=0}{\sum }}M_{i}^{p}\overset{N-1}{\underset{%
j=\left\langle n_{k}\right\rangle +1}{\sum }}\frac{1}{M_{j}^{1-p}}%
+M_{\left\langle n_{k}\right\rangle }^{1-2p}\overset{N-2}{\underset{%
i=\left\langle n_{k}\right\rangle }{\sum }}M_{i}^{p}\overset{N-1}{\underset{%
j=i+1}{\sum }}\frac{1}{M_{j}^{1-p}} \\
&&+c_{p}\overset{N-1}{\underset{i=0}{\sum }}\frac{M_{i}^{p}}{M_{N}^{p}} \\
&\leq &c_{p}M_{\left\langle n_{k}\right\rangle }^{1-2p}M_{\left\langle
n_{k}\right\rangle }^{p}\frac{1}{M_{\left\langle n_{k}\right\rangle }^{1-p}}%
+c_{p}M_{\left\langle n_{k}\right\rangle }^{1-2p}\overset{N-2}{\underset{%
i=\left\langle n_{k}\right\rangle }{\sum }}\frac{1}{M_{i}^{1-2p}}+c_{p}\leq
c_{p}<\infty .
\end{eqnarray*}

The proof of the a) part is complete.

b) Let 
$\left\{ n_{k}:k\geq 0\right\} $ be a sequence of positive numbers,
satisfying condition (\ref{31aaa}). Then 
\begin{equation}
\sup_{k\in \mathbb{N}}\frac{M_{\left\vert n_{k}\right\vert }}{M_{\left\langle
n_{k}\right\rangle }}=\infty .  \label{12h}
\end{equation}

Under condition (\ref{12h}) there exists a sequence $\left\{ \alpha _{k}:%
\text{ }k\geq 0\right\} \subset \left\{ n_{k}:\text{ }k\geq 0\right\} $ such
that $\alpha _{0}\geq 3$ and 
\begin{equation}
\sum_{k=0}^{\infty }\frac{M_{\left\langle \alpha _{k}\right\rangle }^{\left(
1-2p\right) /2}\Phi^{p/2}\left(\alpha_{k}\right)}{M_{\left\vert \alpha _{k}\right\vert }^{\left( 1-2p\right)
/2}}<c<\infty .  \label{12hh}
\end{equation}

Let \qquad 
\begin{equation*}
f^{\left( n\right) }=\sum_{\left\{ k;\text{ }\left\vert \alpha
_{k}\right\vert <n\right\} }\lambda _{k}a_{k},
\end{equation*}%
where 
\begin{equation*}
\lambda _{k}=\frac{\lambda M_{\left\langle \alpha _{k}\right\rangle
}^{\left( 1/p-2\right) /2}\Phi^{1/2}\left(\alpha_{k}\right)}{M_{\left\vert \alpha _{k}\right\vert }^{\left(
1/p-2\right) /2}}
\end{equation*}%
and%
\begin{equation*}
a_{k}=\frac{M_{\left\vert \alpha _{k}\right\vert }^{1/p-1}}{\lambda }\left(
D_{M_{\left\vert \alpha _{k}\right\vert +1}}-D_{M_{\left\vert \alpha
_{k}\right\vert }}\right) .
\end{equation*}

B applying Lemma \ref{lemma1} we can conclude that $f\in H_{p}.$

It is evident that
\begin{equation}
\widehat{f}(j)=\left\{ 
\begin{array}{l}
M_{\left\vert \alpha _{k}\right\vert }^{1/2p}M_{\left\langle \alpha
_{k}\right\rangle }^{\left( 1/p-2\right) /2}\Phi^{1/2}\left(\alpha_{k}\right),\,\,\text{ } \\ 
\text{if \thinspace \thinspace }j\in \left\{ M_{\left\vert \alpha
_{k}\right\vert },...,\text{ ~}M_{\left\vert \alpha _{k}\right\vert
+1}-1\right\} ,\text{ }k=0,1,2..., \\ 
0\text{ },\text{ \thinspace \qquad \thinspace\ \ \ \ \thinspace\ \ \ \ \ }
\\ 
\text{\ if \thinspace \thinspace \thinspace }j\notin
\bigcup\limits_{k=0}^{\infty }\left\{ M_{\left\vert \alpha _{k}\right\vert
},...,\text{ ~}M_{\left\vert \alpha _{k}\right\vert +1}-1\right\} .\text{ }%
\end{array}%
\right.  \label{6aacharp}
\end{equation}

Moreover,%
\begin{equation*}
\frac{\sigma _{_{\alpha _{k}}}f}{\Phi\left(\alpha_{k}\right)}=\frac{1}{\alpha_{k}\Phi\left(\alpha_{k}\right)}\sum_{j=1}^{M_{\left\vert
\alpha _{k}\right\vert }}S_{j}f+\frac{1}{\alpha _{k}\Phi\left(\alpha_{k}\right)}\sum_{j=M_{\left\vert
\alpha _{k}\right\vert }+1}^{\alpha _{k}}S_{j}f:=I+II.
\end{equation*}%
Let $M_{\left\vert \alpha _{k}\right\vert }<j\leq \alpha _{k}.$ Then, by
applying (\ref{6aacharp}) we get that 
\begin{equation}
S_{j}f=S_{M_{\left\vert \alpha _{k}\right\vert }}f+M_{\left\vert \alpha
_{k}\right\vert }^{1/2p}M_{\left\langle \alpha _{k}\right\rangle }^{\left(
1/p-2\right)/2}\Phi^{1/2}\left(\alpha_{k}\right)\left( D_{j}-D_{M_{\left\vert \alpha _{k}\right\vert
}}\right).  \label{8aafn}
\end{equation}

By using (\ref{8aafn}) we can rewrite $II$ as
\begin{eqnarray*}
II &=&\frac{\alpha _{k}-M_{\left\vert \alpha _{k}\right\vert }}{\alpha _{k}\Phi\left(\alpha_{k}\right)}%
S_{M_{\left\vert \alpha _{k}\right\vert }}f+\frac{M_{\left\vert \alpha
_{k}\right\vert }^{1/2p}M_{\left\langle \alpha _{k}\right\rangle }^{\left(
1/p-2\right)/2}}{\alpha _{k}\Phi^{1/2}\left(\alpha_{k}\right)}\sum_{j=M_{\left\vert \alpha _{k}\right\vert
}}^{\alpha _{k}}\left( D_j-D_{M_{\left\vert \alpha _{k}\right\vert
}}\right) \\
&:=&II_{1}+II_{2}.
\end{eqnarray*}

Since (for details see e.g. \cite{BNPT} and \cite{tep9})
\begin{equation*}
\left\Vert S_{M_{\left\vert \alpha _{k}\right\vert }}f\right\Vert
_{weak-L_{p}}\leq c_{p}\left\Vert f\right\Vert_{H_{p}}
\end{equation*}
we obtain that
\begin{eqnarray*}
&&\left\Vert II_{1}\right\Vert _{weak-L_{p}}^{p}\leq\left(\frac{\alpha
	_{k}-M_{\left\vert \alpha _{k}\right\vert }}{\alpha _{k}\Phi\left(\alpha_{k}\right)}\right)
^{p}\left\Vert S_{M_{\left\vert \alpha _{k}\right\vert }}f\right\Vert
_{weak-L_{p}}^{p} \\
&\leq& \left\Vert S_{M_{\left\vert \alpha _{k}\right\vert }}f\right\Vert
_{weak-L_{p}}^{p}\leq c_{p}\left\Vert f\right\Vert _{H_{p}}^{p}<\infty .
\end{eqnarray*}

By using part a) of Theorem \ref{theorem1} we find that
\begin{equation*}
\left\Vert I\right\Vert _{weak-L_{p}}^{p}=\left( \frac{M_{\left\vert \alpha
_{k}\right\vert }}{\alpha _{k}\Phi\left(\alpha_{k}\right)}\right) ^{p}\left\Vert \sigma _{M_{\left\vert
\alpha _{k}\right\vert }}f\right\Vert _{weak-L_{p}}^{p}\leq c_{p}\left\Vert
f\right\Vert _{H_{p}}^{p}<\infty .
\end{equation*}

Let $x\in $ $I_{_{\left\langle \alpha _{k}\right\rangle +1}}^{\left\langle
\alpha _{k}\right\rangle -1,\left\langle \alpha _{k}\right\rangle }.$ Under
condition (\ref{31aaa}) we can conclude that $\left\langle \alpha
_{k}\right\rangle \neq \left\vert \alpha _{k}\right\vert $ and $\left\langle
\alpha _{k}-M_{\left\vert \alpha _{k}\right\vert }\right\rangle
=\left\langle \alpha _{k}\right\rangle .$ Since 
\begin{equation}
D_{j+M_{n}}=D_{M_{n}}+\psi _{M_{n}}D_{j}=D_{M_{n}}+r_{n}D_{j},\text{ when }%
\,\,j<M_{n}  \label{8k}
\end{equation}
if we apply estimate (\ref{9k}) in Lemma \ref{lemma8} for $II_{2}$ we obtain
that 
\begin{eqnarray*}
\left\vert II_{2}\right\vert &=&\frac{M_{\left\vert \alpha _{k}\right\vert
}^{1/2p}M_{\left\langle \alpha _{k}\right\rangle }^{\left( 1/p-2\right) /2}}{
\alpha _{k}\Phi^{1/2}\left(\alpha_{k}\right)}\left\vert \sum_{j=1}^{\alpha _{k}-M_{\left\vert \alpha
_{k}\right\vert }}\left( D_{j+M_{\left\vert \alpha _{k}\right\vert
}}-D_{M_{\left\vert \alpha _{k}\right\vert }}\right) \right\vert \\
&=&\frac{M_{\left\vert \alpha _{k}\right\vert }^{1/2p}M_{\left\langle \alpha
_{k}\right\rangle }^{\left( 1/p-2\right) /2}}{\alpha _{k}\Phi^{1/2}\left(\alpha_{k}\right)}\left\vert \psi
_{M_{\left\vert \alpha _{k}\right\vert }}\sum_{j=1}^{\alpha
_{k}-M_{\left\vert \alpha _{k}\right\vert }}D_{j}\right\vert \\
&\geq &\frac{c_{p}M_{\left\vert \alpha _{k}\right\vert }^{1/2p-1}M_{\left\langle
\alpha _{k}\right\rangle }^{\left( 1/p-2\right) /2}}{\Phi^{1/2}\left(\alpha_{k}\right)}\left( \alpha
_{k}-M_{\left\vert \alpha _{k}\right\vert }\right) \left\vert K_{\alpha
_{k}-M_{\left\vert \alpha _{k}\right\vert }}\right\vert \\
&\geq &\frac{c_{p}M_{\left\vert \alpha _{k}\right\vert }^{1/2p-1}M_{\left\langle
\alpha _{k}\right\rangle }^{\left( 1/p+2\right) /2}}{\Phi^{1/2}\left(\alpha_{k}\right)}.
\end{eqnarray*}

It follows that%
\begin{eqnarray*}
&&\left\Vert II_{2}\right\Vert _{weak-L_{p}}^{p} \\
&\geq &c_{p}\left( \frac{M_{\left\vert \alpha _{k}\right\vert }^{\left(
1/p-2\right) /2}M_{\left\langle \alpha _{k}\right\rangle }^{\left(
1/p+2\right) /2}}{\Phi^{1/2}\left(\alpha_{k}\right)}\right) ^{p}\mu \left\{ x\in G_{m}:\left\vert
IV_{2}\right\vert \geq c_{p}M_{\left\vert \alpha _{k}\right\vert }^{\left(
1/p-2\right) /2}M_{\left\langle \alpha _{k}\right\rangle }^{\left(
1/p+2\right) /2}\right\} \\
&\geq &c_{p}\frac{M_{\left\vert \alpha _{k}\right\vert }^{1/2-p}M_{\left\langle
\alpha _{k}\right\rangle }^{1/2+p}\mu \left\{I_{_{\left\langle \alpha
_{k}\right\rangle +1}}^{\left\langle \alpha _{k}\right\rangle
-1,\left\langle \alpha _{k}\right\rangle }\right\}}{\Phi^{p/2}\left(\alpha_{k}\right)} 
\geq \frac{c_{p}M_{\left\vert \alpha _{k}\right\vert}^{1/2-p}}{M_{\left\langle \alpha
_{k}\right\rangle }^{1/2-p}\Phi^{p/2}\left(\alpha_{k}\right)}.
\end{eqnarray*}
Hence, for large $k$, 
\begin{eqnarray*}
&&\left\Vert \sigma _{\alpha _{k}}f\right\Vert _{weak-L_{p}}^{p} \\
&\geq &\left\Vert II_{2}\right\Vert _{weak-L_{p}}^{p}-\left\Vert
II_{1}\right\Vert _{weak-L_{p}}^{p}-\left\Vert I\right\Vert _{weak-L_{p}}^{p}
\\
&\geq &\frac{1}{2}\left\Vert II_{2}\right\Vert _{weak-L_{p}}^{p}\geq \frac{%
c_{p}M_{\left\vert \alpha _{k}\right\vert }^{1/2-p}}{2M_{\left\langle \alpha
_{k}\right\rangle }^{1/2-p}\Phi^{p/2}\left(\alpha_{k}\right)}\rightarrow \infty ,\text{ as }k\rightarrow
\infty.
\end{eqnarray*}

The proof is complete.
\end{proof}
 
\textbf{Acknowledgment:} We thank the careful referee for some good suggestions which have improved the final version of this paper.


\begin{thebibliography}{99}
\bibitem{AVD} \textit{G. N. Agaev,} \textit{N. Ya. Vilenkin,} \textit{G. M.
Dzafarly} and \textit{A. I. Rubinshtein,} Multiplicative systems of
functions and harmonic analysis on zero-dimensional groups, Baku, Ehim, 1981
(in Russian).

\bibitem{BGG} \textit{I. Blahota,} \textit{G. Gát} and \textit{U. Goginava},
Maximal operators of Fejér means of double Vilenkin-Fourier series, Colloq.
Math. J., 107 (2007), no. 2, 287-296.

\bibitem{BGG2} \textit{I. Blahota,} \textit{G. Gát} and \textit{U. Goginava}%
, Maximal operators of Fejér means of Vilenkin-Fourier series. JIPAM. J.
Inequal. Pure Appl. Math. 7 (2006), 1-7.

\bibitem{bt} \textit{I. Blahota} and \textit{G. Tephnadze,} Strong
convergence theorem for Vilenkin-Fejér means, Publ. Math. Debrecen, 85, 1-2
(2014), 181--196.

\bibitem{BNPT}	I. Blahota, K. Nagy, L. E. Persson, G. Tephnadze, A sharp boundedness result concerning some maximal operators of partial sums with respect to Vilenkin systems, Georgian Math., J., DOI: https://doi.org/10.1515/gmj-2018-0045.

\bibitem{Fu} \textit{N. J. Fujii,} A maximal inequality for $H_{1}$
functions on the generalized Walsh-Paley group, Proc. Amer. Math. Soc. 77
(1979), 111-116.

\bibitem{gat} \textit{G. Gát,} Ces\`{a}ro means of integrable functions with
respect to unbounded Vilenkin systems. J. Approx. Theory 124 (2003), no. 1,
25-43.

\bibitem{GoAMH} \textit{U. Goginava}, Maximal operators of Fejér means of
double Walsh-Fourier series. Acta Math. Hungar. 115 (2007), no. 4, 333-340.

\bibitem{GNCz} \textit{U. Goginava} and \textit{K. Nagy,} On the maximal
operator of Walsh-Kaczmarz-Fejér means, Czechoslovak Math. J., 61 (2011), no. 3, 673-686.

\bibitem{PS} \textit{J. Pál }and \textit{P. Simon, }On a generalization of
the concept of derivative, Acta Math. Acad. Sci. Hungar. 29 (1977), no. 1-2,
155--164.

\bibitem{PT1} \textit{L. E. Persson} and \textit{G. Tephnadze}, A sharp boundedness result concerning some maximal operators of Vilenkin-Fejér means, Mediterr. J. Math., 13, 4 (2016) 1841-1853.

\bibitem{Sc} \textit{F. Schipp,} Certain rearrangements of series in the
Walsh series, Mat. Zametki, 18 (1975), 193-201.

\bibitem{Si1} \textit{P. Simon,} Cesáro summability with respect to
two-parameter Walsh systems, Monatsh. Math., 131, 4 (2000), 321--334.

\bibitem{Si2} \textit{P. Simon,} Investigations with respect to the Vilenkin
system, Ann. Univ. Sci. Budapest. Eötvös Sect. Math., 28 (1985), 87-101.

\bibitem{sm} \textit{B. Smith,} A strong convergence theorem for $%
H_{1}\left( T\right) ,$ Lecture Notes in Math., 995, Springer, Berlin, 1994,
169-173.

\bibitem{tep2} \textit{G. Tephnadze,} On the maximal operator of Vilenkin-Fej%
ér means, Turk. J. Math, 37, (2013), 308-318.

\bibitem{tep3} \textit{G. Tephnadze,} On the maximal operators of
Vilenkin-Fejér means on Hardy spaces, Math. Inequal. Appl., 16, (2013), no.
2, 301-312.

\bibitem{tep8} G. Tephnadze, On the Vilenkin-Fourier coefficients, Georgian Math. J., 20, 1 (2013), 169-177.

\bibitem{tep9}	G. Tephnadze, On the convergence of partial sums with respect to Vilenkin system on the martingale Hardy spaces, J. Contemp. Math. Anal., 53, 5, (2018) 294?306.

\bibitem{Vi} \textit{N. Ya. Vilenkin,} On a class of complete orthonormal
systems, Izv. Akad. Nauk. U.S.S.R., Ser. Mat., 11 (1947), 363-400.

\bibitem{We1} \textit{F. Weisz,} Martingale Hardy spaces and their
applications in Fourier Analysis, Springer, Berlin-Heideiberg-New York, 1994.

\bibitem{We3} \textit{F. Weisz,} Hardy spaces and Ces\`{a}ro means of
two-dimensional Fourier series, Bolyai Soc. Math. Studies, (1996), 353-367.

\bibitem{We2} \textit{F. Weisz,} Cesáro summability of one- and
two-dimensional Walsh-Fourier series, Anal. Math. 22 (1996), no. 3, 229--242.

\bibitem{We4} \textit{F. Weisz,} Weak type inequalities for the Walsh and
bounded Ciesielski systems. Anal. Math. 30 (2004), no. 2, 147-160.
\end{thebibliography}
\end{document}